\documentclass{article}
\pdfoutput=1 

\usepackage{arxiv}
\usepackage{natbib}
\usepackage[utf8]{inputenc} 
\usepackage[T1]{fontenc}    
\usepackage{hyperref}       
\usepackage{url}            
\usepackage{booktabs}       
\usepackage{amsfonts}       
\usepackage{nicefrac}       
\usepackage{microtype}      
\usepackage{graphicx}

\usepackage{doi}

\usepackage[ruled,vlined]{algorithm2e} 

\usepackage{xcolor}
\usepackage{amsmath}
\usepackage{graphicx}
\usepackage{subcaption}

\usepackage{multicol}
\usepackage{tcolorbox}
\usepackage{wrapfig}

\newtheorem{thm}{Theorem}

\newtheorem{assumption}{Assumption}

\newcommand{\bHi}{\bar{G}_i^t}

\newcommand{\tB}{\tilde{B}}

\title{L-DQN: An Asynchronous Limited-Memory Distributed Quasi-Newton Method}

 \author{Bugra Can$^{1}$,\\
 Department of MSIS\\
Rutgers Business School\\
\texttt{bugra.can@rutgers.edu}\\
 \And 
 Saeed Soori$^{2}$,\\
 Department of Computer Science\\
University of Toronto\\
 \texttt{saeed.soori.sh@gmail.com} \\
 \And Maryam Mehri Dehnavi$^{2}$,\\
 Department of Computer Science\\
   University of Toronto \\
   \texttt{mmehride@cs.toronto.edu} \\
 \And Mert G\"urb\"uzbalaban$^{1}$\\
 Department of MSIS \\
   Rutgers Business School\\
   \texttt{mert.gurbuzbalaban@rutgers.edu} 
 \thanks{$^{1}$ Department of Management Sciences and Information Systems, Rutgers Business School, Piscataway, NJ-08854, USA.}
 \thanks{$^{2}$ Department of Computer Sciences, University of Toronto, Toronto, Canada.}
\thanks{*Bugra Can and Mert Gürbüzbalaban acknowledge support from the Office of Naval Research Award Number N00014-21-1-2244, and the grants National Science Foundation (NSF) CCF-1814888, NSF DMS-2053485, NSF DMS-1723085.}
 }
 \date{}

\begin{document}
\maketitle

\begin{abstract}

This work proposes a distributed algorithm for solving empirical risk minimization problems, called L-DQN, under the master/worker communication model. L-DQN is a distributed limited-memory quasi-Newton method that supports asynchronous computations among the worker nodes. 
Our method is efficient both in terms of storage and communication costs, i.e., in every iteration the master node and workers communicate vectors of size $O(d)$, where $d$ is the dimension of the decision variable, and the amount of memory required on each node is $O(md)$, where $m$ is an adjustable parameter. To our knowledge, this is the first distributed quasi-Newton method with provable global linear convergence guarantees in the asynchronous setting where delays between nodes are present. Numerical experiments are provided to illustrate the theory and the practical performance of our method.
	
\end{abstract}

\keywords{Quasi Newton Methods \and Distributed Optimization \and BFGS Methods}

\section{Introduction}

Due to the rapid increase in the size of datasets in the last decade, distributed algorithms that can parallelize the computations to multiple (computational units) nodes connected over a communication network became indispensable \cite{bertsekas1989parallel,recht2011hogwild}. A common communication model in distributed machine learning is the master/worker model in which the master keeps a copy of the global decision variable $x$ and shares it with the workers. 
Each worker operates locally on its own data and then communicates the results to the master to update the decision variable in a synchronous \cite{gurbuzbalaban2017convergence,le2012stochastic,defazio2014saga,defazio2014finito,mairal2015incremental,mokhtari2018surpassing,vanli2018global} or asynchronous fashion \cite{xiao2019dscovr,asaga,arock,bianchi2015coordinate,zhang2014asynchronous,mansoori2017superlinearly,csimcsekli2018asynchronous}. In the synchronous setting, the master waits to receive updates from all workers before updating the decision variable, which can lead to a slow execution if the nodes and/or the network is heterogeneous \cite{kanrar2011performance}. In the asynchronous setting, coordination amongst workers is not needed (or is more relaxed) and the master can proceed with updates without having to wait for slow worker nodes. As a result, asynchronous setting can be more efficient than synchronous in heterogeneous computing environments \cite{wongpanich2020rethinking}. 

In this paper, we consider distributed algorithms for empirical risk minimization, i.e. for solving the finite-sum problem
\vspace{-0.3cm}
\begin{equation}\label{eq:orig_problem}
x^*:= \underset{x\in\mathbb{R}^d}{\mbox{argmin}}\;f(x) := \underset{x\in \mathbb{R}^d}{\mbox{argmin}}\;\frac{1}{n}\sum_{i=1}^{n}f_i(x)
\end{equation} 
where $x\in \mathbb{R}^d$ and $f_i: \mathbb{R}^d\to \mathbb{R}$ is the loss function of node $i\in\{1,...,n\}$. We consider the master/worker communication model with asynchronous computations. With today's distributed computing environments, the cost of communication between nodes is considerably higher than the cost of computation, which leads to sharing matrices of size $O(d^2)$ across nodes to be prohibitively expensive in many machine learning applications. Thus, inspired by prior work \cite{liu1989limited,mokhtari2015global,nash1991numerical,skajaa2010limited,bollapragada2018progressive,berahas2019quasi}, we focus on algorithms that communicate between nodes only vectors of size (at most) $O(d)$. 
There are a number of distributed algorithms for empirical risk minimization that can support asynchronous computations; the most relevant to our work are the recently proposed DAve-RPG \cite{anonymous2018daverpg} and DAve-QN algorithms \cite{SooriDaveQN}. DAve-RPG is a delay tolerant proximal gradient method with linear convergence guarantees that also handles a non-smooth term in the objective. However, it is a first-order method that does not estimate the second-order information of the underlying objective, therefore it can be slow for ill-conditioned problems. DAve-QN is a distributed quasi-Newton method with local superlinear convergence guarantees, however it does not admit global convergence guarantees. Furthermore, it relies on BFGS updates on each node, which requires $O(d^2)$ memory as well as $O(d^2)$ computations for updating the Hessian estimate at each node. For large $d$, this can be slow where DAve-QN looses its edge over first-order approaches \cite{SooriDaveQN}; furthermore its $O(d^2)$ memory requirement can be impractical or prohibitively expensive when $d$ is large, say when $d$ is on the order of ten thousands or hundred thousands.

\textbf{Contributions.} To remedy the shortcomings of the DAve-QN algorithm, we propose L-DQN, a distributed limited-memory quasi-Newton method that requires less memory and computational work per iteration. More specifically, the per iteration and per node storage and computation  of L-DQN are  $O(md)$ and $O(md)$  respectively, where $m$ is a configurable parameter and is the number of vectors stored in the memory at every iteration that contains information about the past gradients and iterates. Because of the reduced storage and computation costs, our proposed algorithm scales well for large datasets, it is communication-efficient as it exchanges vectors of size $O(d)$ at every communication. When the number of nodes is large enough, with an appropriate  stepsize, L-DQN has global linear convergence guarantees for strongly convex objectives, even though the computations are done in an asynchronous manner, as opposed to the DAve-QN method which does not provide global convergence guarantees. In practice, we have also observed that L-DQN works well even if the number of nodes $n$ is not large, for example when $n=2$. To our knowledge, L-DQN is the first distributed quasi-Newton method with provable linear convergence guarantees, even in the presence of asynchronous computations.

\textbf{Related work.} The proposed method can be viewed as an asynchronous distributed variant of the traditional quasi-Newton and limited-memory BFGS methods that have been extensively studied in the optimization community (\cite{goldfarb1970family,broyden1973local,dennis1974characterization,powell1976some}). L-DQN builds on the limited-memory BFGS method \cite{liu1989limited}. Prior work have also investigated incremental gradient (\cite{ig-in,in-mathprog}) and incremental aggregated gradient algorithms (\cite{le2012stochastic,defazio2014saga,defazio2014finito,mairal2015incremental,gurbuzbalaban2017convergence,mokhtari2018surpassing,vanli2018global,mokhtari2018iqn,blatt2007convergent}), which are originally developed for centralized problems. These methods update the global decision variable by processing the gradients of the component functions $f_i$ in a deterministic fashion in a specific (e.g. cyclic) order. They are applicable to our setting in practice, however, these methods do not provide convergence guarantees in asynchronous settings. The Lazily Aggregated Gradient (LAG) \cite{chen2018lag} method, which has a convergence rate similar to batch gradient descent in strongly convex, convex, and nonconvex cases as well as its quantized version \cite{sun2019communication}, is an exception, however, LAG is a first-order method that does not use second-order information. For synchronous settings, the distributed quasi-Newton algorithm proposed by \cite{lee2018distributed} is globally linearly convergent and can handle non-smooth regularization terms; convergence analysis for the algorithm does not exist for asynchronous settings. In this work, we use the star network topology where the nodes follow a master/slave hierarchy. However, there is another setting known as the \emph{decentralized setting} which does not have a master node and communication between the nodes is limited to a given fixed arbitrary network topology  (\cite{nedic2009distributed,mansoori2017superlinearly}). Amongst algorithms for this setting,  \cite{eisen2017decentralized} proposes a linearly convergent decentralized quasi-Newton method and \cite{mansoori2017superlinearly} develops an asynchronous Newton-based approach that has local superlinear convergence guarantees to a neighborhood of the problem 
\eqref{eq:orig_problem}. There are also distributed second-order methods developed for non-convex objectives. Among these, most relevant to our paper are  \cite{csimcsekli2018asynchronous} which proposes a stochastic asynchronous-parallel  L-BFGS method and the DINGO method (\cite{crane2019dingo}) which admits linear convergence guarantees to a local minimum for non-convex objectives that satisfy an invexity property.

\textbf{Notation.} Throughout the paper, we use $\Vert . \Vert$ to denote the matrix 2-norm or the (Euclidean norm) $L_2$ norm depending on the context. The Frobenius norm of a matrix $A\in \mathbb{R}^{n\times m}$ is defined as $\Vert A \Vert^2_F:=\sum_{i=1}^{n}\sum_{j=1}^{m}A_{ij}^2$. The matrix $I_d$ denotes the $d\times d$ identity matrix.{A memory with capacity $m$, denoted as $\mathcal{M}_m$, is a set of tuples $(y,q,\alpha,\beta)$ where $y,q\in \mathbb{R}^d$ and $\alpha,\beta \in \mathbb{R}$; the size of the memory $|\mathcal{M}_m|$ satisfies $|\mathcal{M}_m|\leq m$.}
A function $f: \mathbb{R}^{d}\rightarrow \mathbb{R}$ is called $L$-smooth and $\mu$ strongly convex if for any vector $x,\hat{x} \in \mathbb{R}^{d}$, the Hessian satisfies $\mu \Vert x-\hat{x}\Vert\leq \Vert \nabla^2 f(x)-\nabla^2 f(\hat{x}) \Vert \leq L \Vert x- \hat{x}\Vert$. 
\section{Algorithm}
\subsection{Preliminaries}
\textbf{BFGS algorithm.} In the following, we provide a brief summary of the BFGS algorithm, see \cite{nocedal_num_opt} for more detail. Given a convex smooth function $f:\mathbb{R}^d\rightarrow \mathbb{R}$, the BFGS algorithm consists of  iterations:
\begin{equation*}
x^{t+1}=x^{t}-\eta_t(B^{t+1})^{-1}\nabla f(x^{t}),
\end{equation*}
where $\eta_t$ is a properly chosen stepsize where the matrix $H^{t+1}:=(B^{t+1})^{-1}$ is an estimate of the inverse Hessian matrix at $x^t$ and satisfies the \textit{secant equation}:
\begin{align}\label{eq: secant_eq} 
H^{t+1} y^{t+1}=s^{t+1},
\end{align}
where $s^{t+1}:=x^{t}-x^{t-1}$and $y^{t+1}:=\nabla f(x^{t})-\nabla f(x^{t-1})$ are the differences of the iterates and the gradients respectively. By Taylor's theorem, $y^{t+1}=[\int_0^1 \nabla^2 f(x^{t-1}+\tau (x^{t}-x^{t-1})) d\tau] s^{t+1}$,
therefore for a small enough stepsize $\eta_t$ any matrix $H^{t+1}$ solving the secant equation can be considered as an inverse approximate of Hessian $(\nabla^2 f(x^t))^{-1}$ of the function $f$. In fact, the secant equation \eqref{eq: secant_eq} has infinitely many solutions and quasi-Newton methods differ in how they choose a particular solution. 
BFGS chooses the matrix $H^{t+1}$ according to

\begin{equation*}
H^{t+1}=\left(I-\frac{s^{t+1}(y^{t+1})^{\top}}{(y^{t+1})^\top s^{t+1}}\right)H^t\left(I-\frac{s^{t+1}(y^{t+1})^{\top}}{(y^{t+1})^\top s^{t+1}}\right)+ \frac{s^{t+1}(s^{t+1})^{\top}}{(y^{t+1})^{\top}s^{t+1}}.
\end{equation*}
The corresponding update for $B^{t+1}$ is 
\begin{align}\label{alg: BFGS}
    B^{t+1} &= B^t +U^{t+1} + V^{t+1}, \nonumber \\ U^{t+1}&=\frac{y^{t+1}(y^{t+1})^\top}{(y^{t+1})^\top s^{t+1}},\;\; V^{t+1}= -\frac{B^t s^{t+1}(s^{t+1})^\top B^{t}}{(s^{t+1})^\top B^t s^{t+1}}.
\end{align}
If  function $f$ is strongly convex then $(s^{t+1})^{\top}y^{t+1}>0$ so that the denominator in \eqref{alg: BFGS} cannot be zero. Note that $U^t$ and $V^t$ are both rank-one therefore these updates require $\mathcal{O}(d^2)$ operations. Even though the BFGS algorithm \eqref{alg: BFGS} enjoys local superlinear convergence with an appropriate stepsize, its $\mathcal{O}(d^2)$ memory requirement to store the matrix $B^t$ and $\mathcal{O}(d^2)$ computations required for the updates \eqref{alg: BFGS} may be impractical or prohibitively expensive for machine learning problems when $d$ is large. 

\textbf{Limited-memory BFGS (L-BFGS) algorithm.} Limited-memory BFGS (L-BFGS) requires less memory compared to BFGS algorithm. Instead of storing the whole $B^t$ matrix, L-BFGS stores up to $m$ pairs $\{s^t,y^t\}$ in memory and uses these vectors to approximate the Hessian. The parameter $m$ is adjustable which results in a memory requirement of $\mathcal{O}(md)$. At the start of iteration $t$, we have access to $\{s^j,y^j \}$ for $j=t-m, t-m+1, \dots, t-1$. 
Since the storage is full\footnote{In the beginning of the iterations, when the total number of gradients computed is less than $m$ the storage capacity is not full but the details are omitted for keeping the discussion simpler, see \cite{nocedal_num_opt} for details.}, the oldest pair $\{s^{t-m},y^{t-m} \}$ is replaced by the latest pair $\{ s^t, y^t\}$. The resulting L-BFGS algorithm has the updates: 
\begin{equation*}\label{alg: LBFGS}
x^{t+1}= x^t -\eta_t (\tB^{t+1})^{-1} \nabla f(x^t), 
\end{equation*}
where the matrices $\tB^{t+1}$ are computed 
according to the following formula 

\begin{align*}
    \tilde{B}^{t+1}&= \gamma^{t+1}I_d+\sum_{j=1}^{m} \tilde{U}^{t+2-j}+\tilde{V}^{t+2-j}\nonumber,\\ \tilde{U}^{j}&=\frac{y^{j+1}(y^{j+1})^\top}{(y^{j+1})^\top s^{j+1}}, \;\; \tilde{V}^{j}=-\frac{\tB^{j}s^{j+1}(s^{j+1})^{\top}\tB^{j}}{(s^{j+1})^{\top}\tB^{j}s^{j+1}},
\end{align*}
 where $\gamma^{t+1}$ is a scaling factor. 
We note that L-BFGS requires $\mathcal{O}(md)$ memory which is significantly less compared to $\mathcal{O}(d^2)$ for BFGS for a large $d$. 

\textbf{DAve-QN algorithm.} The DAve-QN algorithm \cite{SooriDaveQN} is an asynchronous quasi-Newton method for solving the optimization problem \eqref{eq:orig_problem} in master/slave communication models. 
Let $x^t$ be the variable that is kept at the master at time $t$ and $z_i^t$ be the local copy that agent $i$ keeps after its last communication with the master. At time $t$, an agent $i_t$ communicates with the master and updates its local estimate $B_{i_t}^t$ for the local Hessian $\nabla^2 f_{i_t}(x^t)$ 
with a BFGS update: 
 \begin{equation} \label{update: LBFGS}
B_{i_t}^{t+1}=B_{i_t}^{t}+\frac{y_{i_t}^{t+1}(y_{i_t}^{t+1})^{\top}}{\alpha^{t+1}_{i_t}}-\frac{q_{i_t}^{t+1}(q_{i_t}^{t+1})^{\top}}{\beta^{t+1}_{i_t}},
 \end{equation} 
 where $q_{i_t}^{t+1}:=B_{i_t}^{t}s_{i_t}^{t+1}$, $y_{i_t}^{t+1}:=\nabla f_i(x^t)-\nabla f_i(z_{i_t}^{t})$, 
  $\alpha^{t+1}_{i_t}:= (y_{i_t}^{t+1})^{\top}s_{i_t}^{t+1}$, and $\beta^{t+1}_{i_t}:=\left(s_{i_t}^{t+1}\right)^{\top}q_{i_t}^{t+1}$
are computed using the local copy $z_{i_t}^t$ and the iterate $x^t$. Let $D_{i_t}^t$ be delay time between information received and send at agent $i_t$ at time t, then agent $i_t$ sends the information $ (B_{i_t}^{t}x^t-B_{i_t}^{t-D_{i_t}^t}z_{i_t}^{t-D_{i_t}^t})$, $y_{i_t}^{t+1}$, $q_{i_t}^{t+1}$, $\alpha_{i_t}^{t+1}$ and $\beta_{i_t}^{t+1}$ to the master after making the update \eqref{update: LBFGS}. Consequently, the master updates the global decision variable with:
$$ 
x^{t}=\left(\sum_{i=1}^{n}B_{i}^t\right)^{-1}\left[\sum_{i=1}^{n}B_i^tz_i^t-\nabla f_i(z_{i}^t)\right].
$$
In the next section, we introduce the L-DQN method which is a limited-memory version of the DAve-QN algorithm. L-DQN will allow us to improve performance for large dimensional problems.
The basic idea is that each agent stores $m$-many tuples $\{y_i^j,q_i^j,\alpha_i^j,\beta_i^j\}$
requiring $O(md)$ memory instead of storing the $d\times d$ matrix $B_{i}^{t}$ and carries out L-BFGS-type updates \eqref{update: LBFGS} to compute the Hessian estimate $\nabla^{2} f_{i}(x^{t})$.

\subsection{A Limited-Memory Distributed Quasi-Newton Method (L-DQN)}
\begin{figure}
    \centering
    \includegraphics[width=0.6\linewidth]{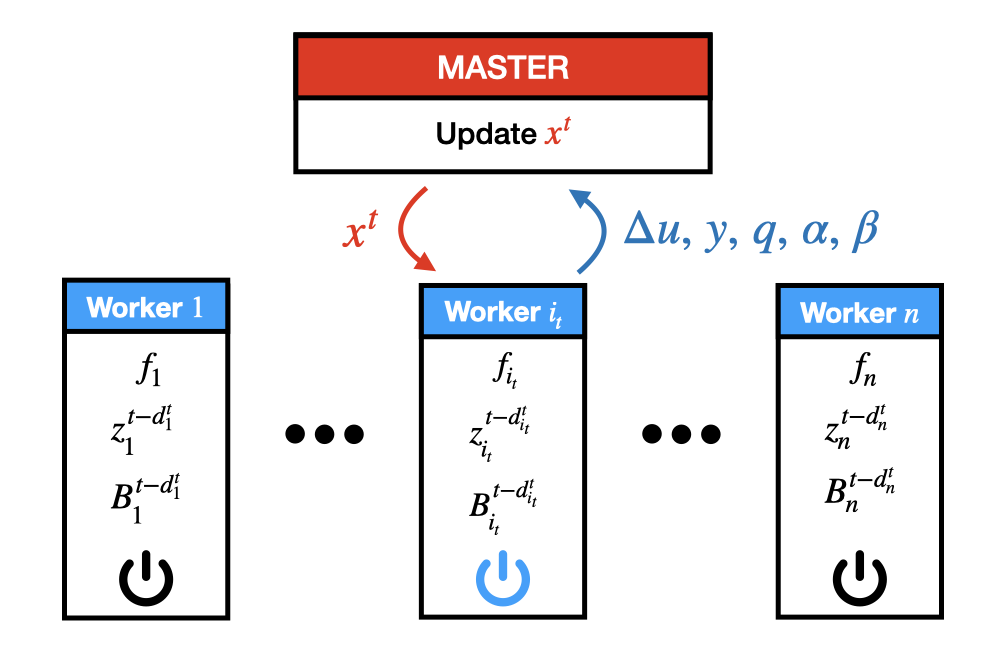}
    \caption{Asynchronous communication scheme used by proposed algorithm.}
    \label{fig:L-DQN}

\end{figure}

In this section, we introduce the L-DQN algorithm on a master/slave communication setting that consists of $n$ workers that are connected to one master with a star topology (see Figure \ref{fig:L-DQN}). Let $d_i^t$ be the delay in communication at time $t$ with the $i$-th worker and the master and $D_i^t$ denote the (penultimate) double delay in communication, i.e. the last exchange between the master and the worker $i$ was at time $t-d_i^t$, and before that the communication took place at $t-D_i^t$ where $D_i^t=d_i^t+d_i^{t-d_i^t-1}+1$. For example, if the node $i$ communicated with master at times $t=1, t= 5$ and $t=7$, we have $d_i^4=3, d_i^5 = 0, d_i^6=1, d_i^7 = 0, d_i^8 = 1$ and $D_i^6 = 6$, $D_i^7 = 2$ and $D_i^8 = 3$.

Let us introduce the historical time $T_i(t)=t-d_i^t$ with the convention $T_i^0(t)=t$ and $T_i^n(t)=T_i(T^{n-1}_i(t))$. We introduce the notation $\tilde{q}_{i_t}^{t+1}:=\tB_{i_t}^{t}s_{i_t}^{t+1}$, $\tilde{\beta}^{t+1}_{i_t}:=\left(s_{i_t}^{t+1}\right)^{\top}\tilde{q}_{i_t}^{t+1}$, and explain the L-DQN updates on worker and master in detail: 

\textbf{Worker Updates:} Each agent $i$ keeps $m$-many tuples $\{y_i^j,\tilde{q}_i^j,\alpha_i^j,\tilde{\beta}_i^j\}$ at their local memory $\mathcal{M}_{m}(i,t)$ at time $t$ and 
at the end of the $m$-th iteration, the worker replaces the oldest tuple  $\{y_i^{T_{i}^m(t)},\tilde{q}_{i_t}^{T_{i}^{m}(t)},\alpha_i^{T_{i}^m(t)},\tilde{\beta}_i^{T_{i}^{m}(t)} \}$ with the new one $\{y_i^{T_{i}(t)},\tilde{q}_{i}^{T_{i}(t)},\alpha_i^{T_{i}(t)},\tilde{\beta}_i^{T_{i}(t)}\}$. Suppose master communicates with worker $i_t$ at the moment $t-d_{i_t}^t$ and sends the copy $x^{t-d_{i_t}^t}$; then upon receiving $x^{t-d_{i_t}^t}$, the worker  $i_t$ computes $\tB_{i_t}^{t+1}x^{t-d_{i_t}^t}$ 
{where $\tB_{i_t}^{t+1}$ is computed according to} 
\begin{small}
\begin{align}\label{update: Bi_LDAveQN}
    \tilde{B}^{t+1}_{i_t}= \gamma^{t+1}_{i_t}I_d+\sum_{j=1}^{m} \tilde{U}^{T_{i_t}^j(t+1)}_{i_t}+\tilde{V}^{T_{i_t}^j(t+1)}_{i_t}\nonumber,\\ \tilde{U}^{j}_{i_t}=\frac{y^{j}_{i_t}(y^{j}_{i_t})^\top}{(y_{i_t}^{j})^\top s^{j}_{i_t}}, \;\; \tilde{V}^{j}_{i_t}=-\frac{\tilde{q}_{i_t}^{j}\left(\tilde{q}_{i_t}^{j}\right)^\top}{(s_{i_t}^{j})^{\top}\tB^{j}s_{i_t}^{j}},
\end{align}
\end{small}
and the scaling factor is chosen as
$\gamma_{i_t}^{t+1}=\frac{\Vert y_{i_t}^{t+1}\Vert^2}{(y_{i_t}^{t+1})^\top s_{i_t}^{t+1}}$.

A number of choices for $\gamma_{i}^t$ are proposed in the literature \cite{nocedal_num_opt}. 
$\gamma_i^t$ given above (which is also considered at \cite{mokhtari2015global}) is an estimate for the largest eigenvalue of Hessian $\nabla^2 f_i(x^{t-d_i^t})$ and works well in practice, therefore our algorithm analysis is based on given $\gamma_i^t$. However, our analysis on the linear convergence of L-DQN can be extended to different choice of $\gamma_i^t$'s as well. 

{
\SetAlFnt{\footnotesize}
\begin{algorithm}[h]
\SetKwFunction{LBFGS}{LBFGS}
\textbf{Function:}\,$u$\,=\,\LBFGS{$\gamma,\mathcal{M}_m,x$}\\
\SetAlgoLined
Set $u=\gamma x$\\
\For{i=1,..,m}{
Retract $y^i,\tilde{q}^i,\alpha^i,\tilde{\beta}^i$ from $\mathcal{M}_{m}$\\
Set $c_1=\frac{(y^i)^{\top}x}{\alpha^i}$ and $c_2=\frac{(\tilde{q}^i)^\top x}{\tilde{\beta}^i}$\\
$u=u+c_1y^i-c_2\tilde{q}^i$
}
\textbf{Return} $u$.
\caption{Compute $u$ given memory $\mathcal{M}_{m} := \{y^i, \tilde{q}^i, \alpha^i, \tilde{\beta}^i\}_{i=1}^m$}
\label{alg: L-Bi Updates }
\end{algorithm}

Worker $i_t$ calls Algorithm \ref{alg: L-Bi Updates } to perform the update \eqref{update: Bi_LDAveQN} locally based on its memory $\mathcal{M}_m(i,t)$. Then, the worker sends
$\Delta u_{i_t}^{t+1}:=\tB_{i_t}^{t+1}x^{t}-\tB_{i_t}^{t-d_{i_t}^t}z_{i_t}^{t-d_{i_t}^t}$, $y_{i_t}^{t-d_{i_t}^t}$, $\tilde{q}_{i_t}^{t-d_{i_t}^t}$, $\alpha_{i_t}^{t-d_{i_t}^t}$ and  $\tilde{\beta}_{i_t}^{t-d_{i_t}^t} 
$ to the master.

\textbf{Master Updates:} Following its communication with the worker, the master receives the vectors $\Delta u_{i_t}$, $y_{i_t}$, $\tilde{q}_{i_t}$, the scalars  $\alpha_{i_t}$, $\tilde{\beta}_{i_t}$ and computes 
\begin{small}
\begin{equation} \label{update: Iterate_LDAveQN_Master}
x^{t+1}=\left(\tB^t \right)^{-1} \left[ \sum_{i=1}^{n}\tB_i^tz_i^t-\eta_t\sum_{i=1}^{n}\nabla f_i(z_i^t) \right],
\end{equation}
\end{small}
where $\tB^t := \sum_{i=1}^{n} \tB_i^t = \sum_{i=1}^{n}\tB_i^{t-d_i^t}$ and stepsize $\eta_t$ determined by the master. Soori et al. have shown in \cite{SooriDaveQN} that 
the computation of $\tB^t$ and $(\tB^{t})^{-1}$ can be done at master locally by using only vectors send by workers. In particular, if we define $u^t:= \sum_{i=1}^{n}\tB_i^tz_i^t=\sum_{i=1}^{n}\tB_i^{t-d_i^t}z_i^{t-d_i^t}$ and
$g^t:=\sum_{i=1}^{n}\nabla f_i(z_i^t)=\sum_{i=1}^{n}\nabla f_i(z_i^{t-d_i^t})$,
then the updates at the master follow the below rules: 
    $\tB^{t+1}=\tB^{t}+(\tB_{i_t}^{t}-\tB_{i_t}^{t-d_i^t})$,
    $u^{t+1}= u^t + \left(\tB^{t+1}_{i_t}x^{t-d_{i_t}^t}-\tB_{i_t}^{t-d_{i_t}^t}x^{t-D_{i_t}^t} \right)$, and $
    g^{t+1}
    =g^t+\left( \nabla f_{i_t}(x^{t-d_{i_t}^t})-\nabla f_{i_t}(x^{t-D_{i_t}^t}) \right)$.
Hence the master only requires $\tB_{i_t}^{t+1}$ and $\nabla f_{i_t}(z_{i_t}^{t+1})=\nabla f_{i_t}(x^{t-d_{i_t}^t})$ to proceed to $t+1$. Let 
\begin{small}
\begin{equation}\label{master: update1} 
U^{t+1}:= (\tB^t)^{-1}-\frac{(\tB^t)^{-1}y_{i_t}^{t+1}(y_{i_t}^{t+1})^{\top}(\tB^t)^{-1}}{(y_{i_t}^{t+1})^\top s_{i_t}^{t+1}+(y_{i_t}^{t+1})^{\top}(\tB^t)^{-1}y_{i_t}^{t+1}},
\end{equation}
\end{small}
then Sherman-Morrison-Woodbury formula implies 

\begin{footnotesize}
\begin{equation} \label{master: update2}
\left(\tB^{t+1}\right)^{-1}= U^{t+1}+\frac{U^{t+1}(\tB_{i_t}^{t-d_{i_t}^t}s_{i_t}^{t+1})(\tB_{i_t}^{t-d_{i_t}^t}s_{i_t}^{t+1})^{\top}U^{t+1}}{(s_{i_t}^{t+1})^{\top}\tB_{i_t}^{t-d_{i_t}^t}s_{i_t}^{t+1}-(\tB_{i_t}^{t-d_i^t}s_{i_t}^{t+1})^{\top}U^{t+1}(\tB_{i_t}^{t-d_i^t}s_{i_t}^{t+1})}.
\end{equation} 
\end{footnotesize}
Thus, if the master already has $(\tB^t)^{-1}$, then $(\tB^{t+1})^{-1}$ is computed using the vectors $y_{i_t}^{t+1}$ and $\tilde{w}^{t+1}:=U^{t+1}\tilde{q}^{t+1}_{i_t}$.
\setlength{\algomargin}{0em}
\SetAlFnt{\footnotesize}
\begin{algorithm*}[ht]

\caption{L-DQN}
\label{Alg: LDQN}
\begin{multicols}{2}
\begin{tcolorbox}[colframe=blue!80,colback=white, title= Worker i:]
\SetKwFunction{LBFGS}{LBFGS}

\textbf{Initialize} $x_i=x^0$,$y_i=x^0$ ,$\tB_i$, $u_{-1}=0$ and memory $\mathcal{M}_{m}=\{ \}$ with capacity $m$.\\
\While{ not interrupted by master}{
  {\color{red} Receive x from master}\\
  $s_i=x-z_i $, $y_i=\nabla f_i(x)-\nabla f_i(z_i)$, $\gamma_i=\frac{y_i^\top y_i}{s_i^\top y_i}$\\
  Compute $\tilde{q}_i$=\LBFGS{$\gamma_i,\mathcal{M}_m,s_i$}\\
  $\alpha_i=y_i^{\top}s_i$, \,$\tilde{\beta}_i=s_i^{\top}\tilde{q}_i$\\
  \uIf {$|\mathcal{M}_{m}| < m$}{
  Add $\{y_i,\tilde{q}_i,\alpha_i,\tilde{\beta}_i \}$ to $\mathcal{M}_{m}$
  }\ElseIf {$|\mathcal{M}_{m}|=m$}{
  Replace the oldest tuple with $(y_i,\tilde{q}_i,\alpha_i,\tilde{\beta}_i)$ at $\mathcal{M}_m$
  }
  Compute $u$=\LBFGS{$\gamma_i,\mathcal{M}_m,x$}\\
  $\Delta u=u-u_{-1}$\\ 
  $u_{-1}=u$, $z_i=x$\\
  {\color{blue!80} Send $\Delta u, y_i,\tilde{q}_i,\alpha_i,\tilde{\beta}_i$ to the master}
}

\end{tcolorbox}
\columnbreak
\begin{tcolorbox}[colframe=red!80,colback=white ,title= Master:]
\begin{small}

\textbf{Initialize} $x$, $\eta$, $\tB_i$, $g= \sum_{i=1}^{n}\nabla f_i(x)$, $\tB^{-1}=(\sum_{i=1}^{n}\tB_i)^{-1}$, $u = \sum_{i=1}^{n}\tB_ix$.\\
\For{$t=1,...,T$}{
    \textbf{If} a worker sends an update\\
    {\color{blue} Receive $\Delta u, y,\tilde{q},\alpha,\tilde{\beta}$ from worker}\\
    $u=u+\Delta u,\; g=g+y,\; v=\tB^{-1}y$\\
    $U=\tB^{-1}-\frac{vv^\top}{\alpha+v^\top y}$\\
    $w= U\tilde{q},\; \tB^{-1}=U+\frac{ww^\top}{\tilde{\beta}-\tilde{q}^\top w}$\\
    $x=\tB^{-1}(u-\eta g)$\\
    {\color{red!80} Send x to the worker in return}
    }
Interrupt all the workers.\\
\textbf{Output} $x^T$
\end{small}

\end{tcolorbox}
\end{multicols}

\end{algorithm*}

}
{
The steps for the master  and worker nodes are provided in Algorithm  \ref{Alg: LDQN}. After receiving $x^t$ from the master, worker $i$ computes its estimate $\tB_i^{t+1}$ using the vector $x^{t-d_i^t}$ received from the master, then updates its memory $\mathcal{M}_m(i,t)$ and sends the vectors $\Delta u_i, y_i ,\tilde{q}_i$ together with  the scalars $\alpha_i, \tilde{\beta}_i$ back to the master. Based on \eqref{master: update1} and \eqref{master: update2}, the master computes $x^{t+1}$ using the vectors received from worker $i$.
}We define the \textit{epochs} $\{E_m\}_{m\in \mathbb{N}^{+}}$ recursively as follows: We set $E_1=0$ and define
$
E_{m+1}= \min \{ t: t-D_i^t \geq E_m \;\mbox{for all}\; i=1,...,n \}.
$
In other words, $E_{m+1}$ is the first time $t$ such that each machine makes at least 2 updates on the interval $[E_m,t]$. Epochs as defined above satisfy the properties: 
\begin{itemize}
    \item For any $t \in [E_{m+1},E_{m+2})$ and any $i=1,2,..,N$ one has $t-D_i^t \in [E_{m},t)$
    \item If delays are uniformly bounded, i.e. there exists a constant $d_i^t\leq d$ for all $i$ and $t$, then for all $m$ we have $E_{m+1}-E_m \leq D:= 2d+1$ and $E_m \leq Dm$.
    \item If we define average delays as $\bar{d^t}:= \frac{1}{N}\sum_{i=1}^{n}d_i^t$, then $\bar{d^t}\geq (n-1)/2$. Moreover, assuming that $\bar{d^t}\leq (n-1)/2+\bar{d}$ for all t, we get $E_{m}\leq 4n(d+1)m$.
\end{itemize}
Notice that convergence to optimum $x^*$ is not possible without visiting every function $f_i$, so measuring performance using epochs where every node has communicated with the master at least once is more convenient than the number of communications, $t$, for comparison.

\section{Convergence Analysis}
In this section, we study theoretical results for linear convergence of L-DQN algorithm with a constant stepsize $\eta_t=\eta$. Firstly, we assumed that the functions $f_i$'s and the matrices $\tB_i^t$'s satisfy the following conditions:
\begin{assumption}\label{Assump 1} The component functions $f_i$ are $L$- smooth and $\mu$-strongly convex, i.e. there exist positive constants $0<\mu<L$ such that, for all $i$ and $x,\hat{x}\in \mathbb{R}^p$,
$$
\mu \Vert x-\hat{x}\Vert^2 \leq (\nabla f_i(x)-\nabla f_i(\hat{x}))^\top(x-\hat{x})\leq L\Vert x-\hat{x}\Vert^2.
$$
\end{assumption} 
\begin{assumption}\label{Assump 2} 
There exist constants $0 < \epsilon_d <\epsilon_u$ such that the following bounds are satisfied for all $i=1,...,n$ and $x\in\mathbb{R}^d$ at any $t>0$:
\begin{equation}\label{Assump: Bounds} 
\epsilon_d I_d \preceq (\tB_i^t)^{-1/2}\nabla^2f_i(x^t)(\tB_i^{t})^{-1/2} \preceq \epsilon_u I_d.
\end{equation}
\end{assumption}
 
Assumption \ref{Assump 2} says that $\tB_i^t$ approximates the Hessian $\nabla^2 f_i(x^t)$ up to a constant factor.For example, if the objective is a quadratic function of the form $f_i(x)=\frac{1}{2}(x-x_*)^\top Q_i (x-x_*)$ and $\frac{1}{1+c} Q_i \preceq \tB_i^t \preceq  (1+c) Q_i$ for some constant $c>0$ then we would have $\epsilon_d=\frac{1}{1+c}, \epsilon_u = (1+c)$ and the ratio $\epsilon=\epsilon_u/\epsilon_d$
satisfies $\epsilon = (1+c)^2\geq 1$. In fact, this ratio can be thought as a measure of the accuracy of the Hessian approximation. In the special case when the Hessian approximations are accurate (when $\tB_i^t =  Q_i$), we have $c=0$ and $\epsilon=1$. Otherwise, we have $\epsilon>1$.

In particular, if the eigenvalues of $\tB_i^t$ stay in the interval $[\lambda_d,\lambda_u]$, then Assumption \eqref{Assump 2} holds with $\epsilon_d=\frac{\mu}{\lambda_u}$ and $\epsilon_{u}=\frac{L}{\lambda_d}$.
We note that in the literature, there exist estimates for $\lambda_u$ and $\lambda_d$ \cite{nocedal_num_opt,mokhtari2015global}.
For example, it is known that if we choose $\gamma_i^t =\frac{\Vert y_i^{t}\Vert^2}{(y_{i}^{t})^{\top}s_i^{t}}$ then we can take $\lambda_u=(m+d)L$, and  $ \lambda_d=\frac{(\mu)^{m+d}}{((m +d)L)^{m+d-1}}$ for memory/storage capacity $m$ (see \cite{mokhtari2015global} for details). A shortcoming of these existing bounds \cite{efficientlycomputeeigen,LargeScaleTRS} is that they are not tight, i.e. with an increasing memory capacity $m$, the bounds get worse.
In our experiments, we have also observed in real datasets that Assumption \ref{Assump 2} holds where we estimated the constants $\epsilon_d$ and $\epsilon_u$ (see the supplementary file for details). These numerical results show that Assumption \ref{Assump 2} is a reasonable assumption to make in practice for analyzing L-DQN methods.

Before we provide a convergence result for L-DQN, we observe that the iterates $\{x^t\}_{t \in \mathbb{N}^{+}}$ of L-DQN provided in \eqref{update: Iterate_LDAveQN_Master} satisfy the property 

\begin{footnotesize}
\begin{equation}\label{eq: Iterates}
x^{t+1}-x^* =\left(\tB^t \right)^{-1} \left[\sum_{i=1}^{n} \tB_i^t (z_i^t-x^*)-\eta_t \sum_{i=1}^{n}\nabla f_i(z_i^t)-\nabla f_i(x^*) \right].
\end{equation} 
\end{footnotesize}
The next theorem uses bounds  \eqref{Assump: Bounds} together with equality \eqref{eq: Iterates} to find the condition on fixed step size $\eta := \eta_t$ such that the L-DQN algorithm is linearly convergent on epochs $[E_m,E_{m+1})$ for $m\in \mathbb{N}^{+}$. The proof can be found in the appendix.
  
\begin{thm}\label{thm: L_DAveQN Convergence}
Suppose Assumptions \ref{Assump 1}-\ref{Assump 2} hold and accuracy $\epsilon:= \epsilon_u/\epsilon_d$ 
satisfies 
\begin{footnotesize}
\begin{equation}\label{cond: epsilon_wo_n}
\epsilon < \frac{1}{2}\left[1+\frac{1}{\kappa}+\sqrt{\left(1+\frac{1}{\kappa} \right)^2+\frac{4}{\kappa}} \right],
\end{equation}
\end{footnotesize}
where $\kappa=\frac{L}{\mu}$ is the condition number of $f$. Let stepsize $\eta \in (\frac{1}{\epsilon_d}[1-\frac{1}{\epsilon \kappa}], \frac{2}{\epsilon_d+\epsilon_u})$, then for each $t \in [E_{m-1}, E_{m})$, the iterates $x^t$ generated by L-DQN algorithm \eqref{Alg: LDQN} converge linearly on epochs $[E_m,E_{m+1})$, i.e. there exists $\rho<1$ such that $\Vert x^{t+1}-x_*\Vert\leq \rho^{m} \Vert x^{0}-x_*\Vert$ for all $t \in [E_m,E_{m+1})$.
\end{thm}

\setlength{\dbltextfloatsep}{0pt}
Theorem \ref{thm: L_DAveQN Convergence} says that if the Hessian approximations are good enough, then $\varepsilon$ is small enough and L-DQN will admit a linear convergence rate. 
Even though the condition \eqref{cond: epsilon_wo_n} seems conservative on the accuracy of the Hessian estimates, to our knowledge, there exists no other linear convergence result that supports global convergence of asynchronous distributed second-order methods for distributed empirical risk minimization. Also, on real datasets, we observe that L-DQN algorithm performs well even though limited memory updates fail to satisfy the condition \eqref{cond: epsilon_wo_n} on the accuracy. 

\begin{figure*}[h!]
    \centering
    \begin{subfigure}[b]{0.32\textwidth}
        \includegraphics[width=\textwidth]{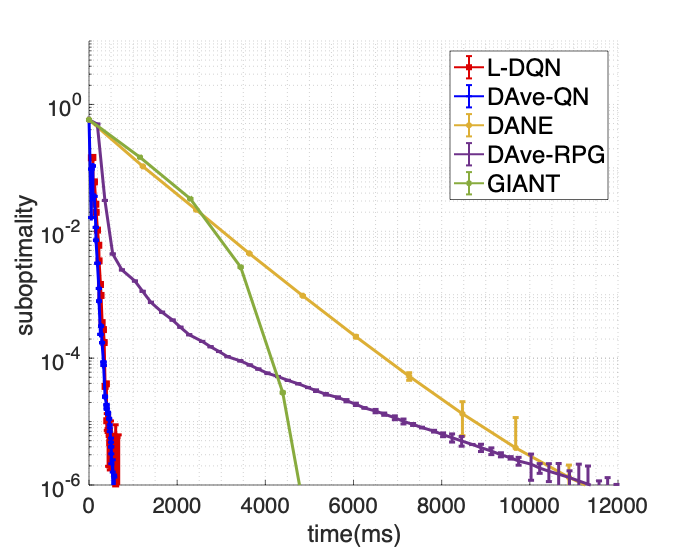}
         \caption{covtype$(54,\approx2.9M)$}
        \label{fig:cov}
    \end{subfigure}
    \begin{subfigure}[b]{0.32\textwidth}
        \includegraphics[width=\textwidth]{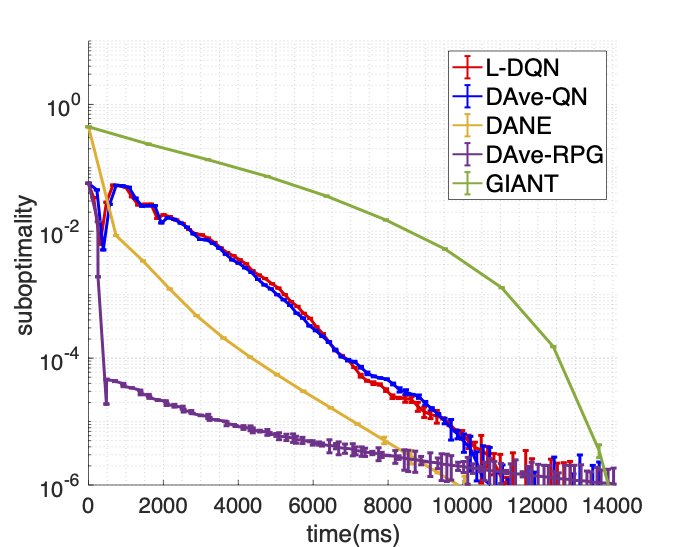}
         \caption{SVHN$(3072,\approx73K)$}
        \label{fig:susy}
    \end{subfigure}
    \begin{subfigure}[b]{0.32\textwidth}
        \includegraphics[width=\textwidth]{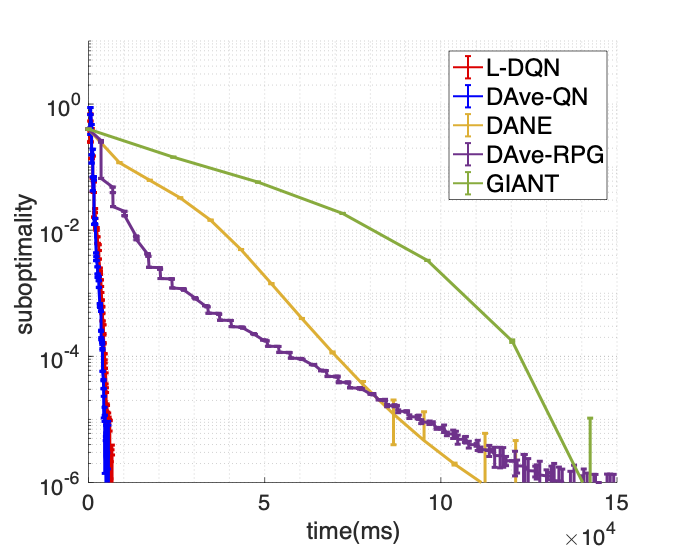}
         \caption{mnist8m$(784,8.1M)$}
        \label{fig:mnist8m}
    \end{subfigure}
    \caption{\small Expected suboptimality versus time. The first and second numbers adjacent to the dataset names are variables $d$ and $n$ respectively. }\label{fig:suboptimality}

\end{figure*}

\section{Numerical Experiments}\label{sec: Experiments}

We tested our algorithm on the multi-class logistic regression problem with $L_2$ regularization where the objective is

\begin{small}
\begin{equation}
    f(x) = \frac{1}{n} \sum_{i=1}^n  \log(1+\exp(-b_i a_i^T x)) +\frac{\lambda}{2}\|x\|^2.
\end{equation}
\end{small}
and $\lambda>0$ is the regularization parameter, $a_i\in \mathbb{R}^d$ is a feature vector, and $b_i$ is the corresponding label. We worked with five datasets (\verb+SVHN+, \verb+mnist8m+, \verb+covtype+, \verb+cifar10+, \verb+rcv1+) from the LIBSVM repository \cite{chang2011libsvm} where the \verb+covtype+ dataset is expanded based on the approach in \cite{wang2018giant} for large-scale experiments.

We compare L-DQN with the following other recent distributed optimization algorithms:
\begin{itemize}
\item \textbf{DAve-QN} \cite{SooriDaveQN}: An asynchronous distributed quasi-Newton method. 
\item \textbf{Distributed Average Repeated Proximal Gradient (DAve-RPG)}\cite{mishchenko2018delay}: A first-order asynchronous method that performs better compared to incremental aggregated gradient \cite{gurbuzbalaban2017convergence} and synchronous proximal gradient methods.
\item \textbf{Globally Improved Approximate Newton (GIANT)}\cite{giant2018mahoney}: A synchronous communication-efficient approximate Newton method, for which the numerical experiments in \cite{giant2018mahoney} demonstrate it outperforms DANE \cite{shamir2014dane} and the Accelerated Gradient Descent \cite{nesterov2013introductory}. 
\item \textbf{Distributed Approximate Newton (DANE)}\cite{shamir2014dane}:  A well-known second-order method that requires sychronization step among all workers. 


\end{itemize}

All experiments are conducted on the XSEDE Comet resources \cite{towns2014xsede} with 24 workers (on Intel Xeon E5-2680v3 2.5 GHz architectures) and with 120 GB Random Access Memory (RAM.) For L-DQN, DAve-RPG, and DAve-QN, we use 17 processes where one as a master and the 16 processes dedicated as workers. DANE and GIANT do not have a master, thus, we use 16 processes are workers with no master. All datasets are normalized to [0,1] and randomly distributed so the load is roughly balanced among workers. We use Intel MKL 11.1.2 and MVAPICH2/2.1 for BLAS (sparse/dense) operations and MPI programming compiled with mpicc 14.0.2 for optimized  communication.
Each experiment is repeated five times and the average and standard deviation is reported as error bars in our results. 

\textbf{Parameters}: The recommended parameters for each method is used.  $\lambda$ is tuned to ensure convergence for all methods. We use $\lambda=1$ for \verb+mnist8m+, and $\lambda=0.1$ for \verb+SVHN+, \verb+cifar10+ and \verb+covtype+. Other choices of $\lambda$ show similar performances. For DANE, SVRG \cite{johnson2013accelerating} is used as a local solver; parameters are selected based on experiments in \cite{shamir2014communication}. DANE has two parameters $\eta$ and $\mu$ which are set to 1 and $3\lambda$ respectively based on the recommendation of the authors in \cite{shamir2014communication}. For DAve-RPG, the number of passes on local data is set to 5 ($p=5$) and its stepsize is selected using a standard backtracking line algorithm \cite{schmidt2015non}. For L-DQN, the memory capacity is set as $m=20$ for \verb+covtype+, \verb+mnist8m+ and \verb+cifar10+ where stepsize is $\eta=0.8$,$\eta=0.8$ and $\eta=0.6$ respectively. On \verb+SVHN+, the parameters of L-DQN are chosen as $m=25$ and $\eta=0.9$.

\begin{figure}[ht!]
    \centering
    \includegraphics[width=0.49\linewidth]{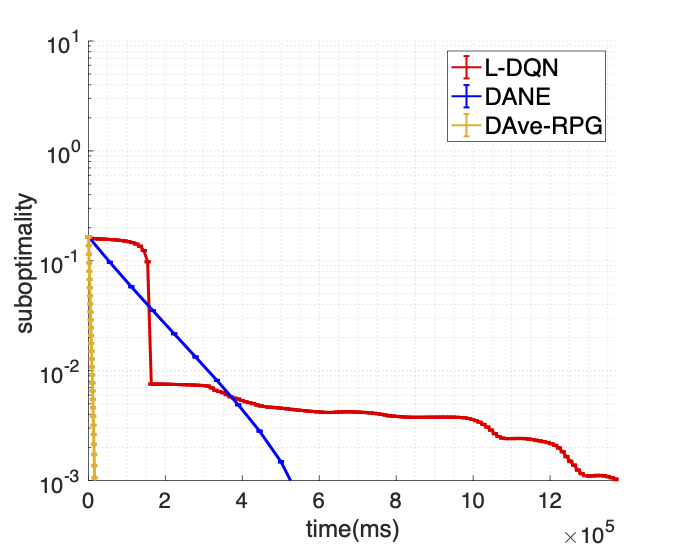}
    \includegraphics[width=0.49\linewidth]{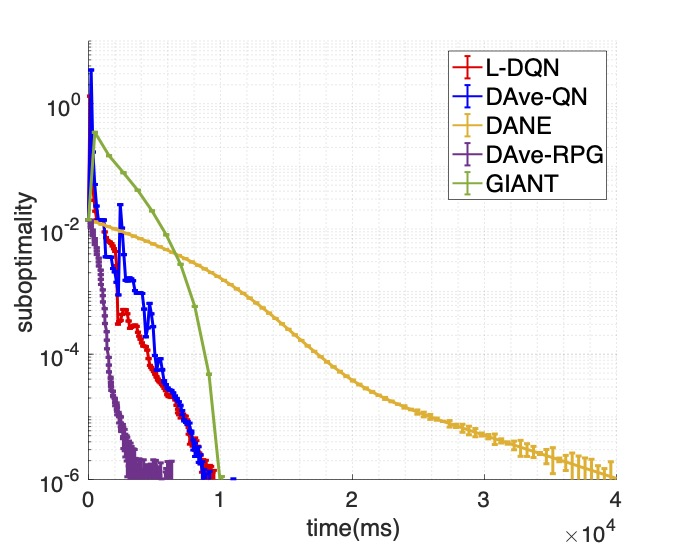}
    \caption{\small Expected suboptimality versus time on rcv1 (left) and cifar10 (right)
    }\label{fig:suboptimality-extra}
\end{figure}

\autoref{fig:suboptimality} shows the average suboptimality versus time for the datasets \verb+mnist8m+, \verb+SVHN+ and \verb+covtype+. We observe that L-DQN converges with a similar rate compared to DAve-QN while it uses less memory. For larger datasets (such as the \verb+rcv1+ with $d=47,000$ and $n=697,000$ at Figure \ref{fig:suboptimality-extra}), DAve-QN was not able to run due to its memory requirement whereas the other methods run successfully. DAve-RPG demonstrates good performance at the beginning for  \verb+SVHN+ compared to other methods due to its cheaper iteration complexity. However, L-DQN becomes faster eventually and outperforms DAve-RPG.

The right panel of \autoref{fig:suboptimality-extra} shows the suboptimality versus time for the dataset \verb+cifar10+ where we choose the parameter
$\lambda=0.1$,  $m=20$ and $\eta=0.6$ for \verb+cifar10+. DAve-RPG is the fastest on this dataset whereas L-DQN is competitive with DAve-QN with less memory requirements. We conclude that when the underlying optimization problem is ill-conditioned (such as the case of \verb+mnist8m+ dataset), L-DQN improves performance with respect to other methods while being scalable to large datasets. In case of less ill-conditioned problems (such as \verb+SVHN+ and \verb+cifar10+), first-order methods such as DAve-RPG are efficient where second-order methods may not be necessary.

\begin{figure}[ht!]
    \centering
    \includegraphics[width=1.05\linewidth]{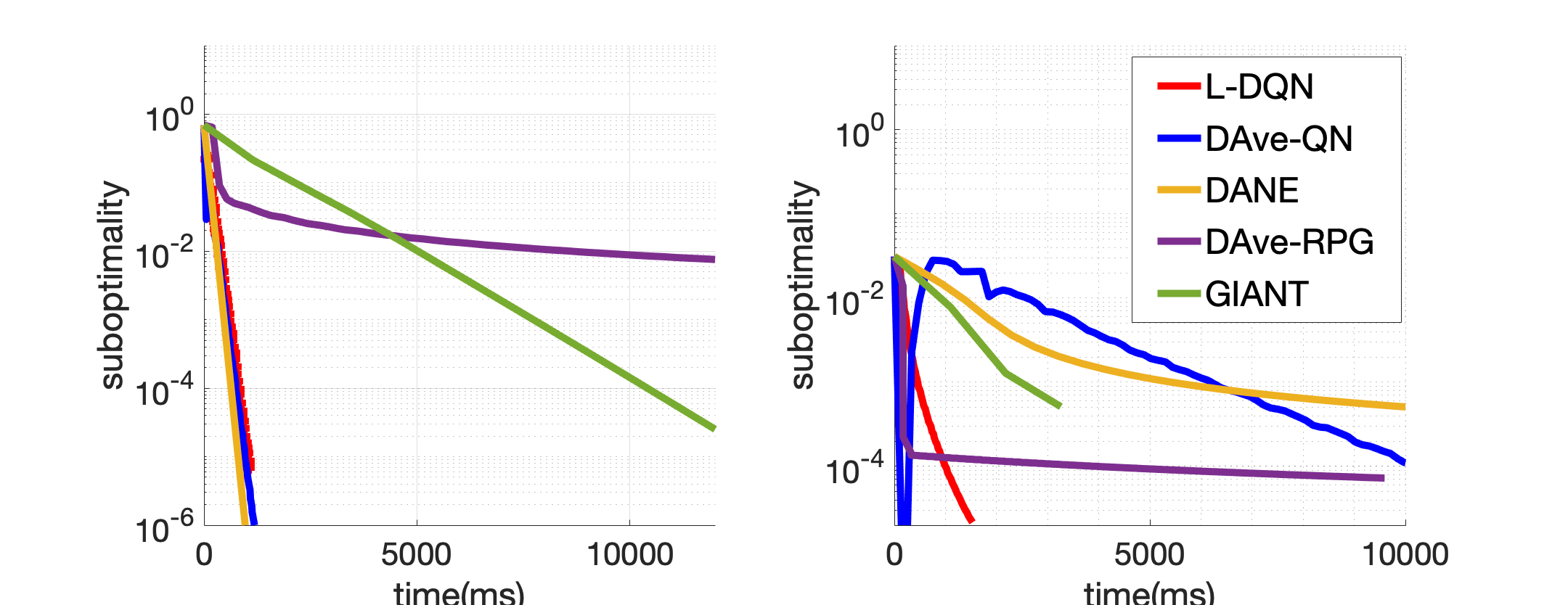}
    \caption{\small Suboptimality comparison without strong convexity assumption on datasets \texttt{covtype}(left) and \texttt{cifar10}(right).}
    \label{fig:LDQN_vs_DAveQN}

\end{figure}
Figure \ref{fig:LDQN_vs_DAveQN} exhibits the suboptimality results of the algorithms on \verb+cifar10+ and \verb+covtype+ without regularization parameter which makes the problems more ill-conditioned. Due to its less memory requirement, we can see that the performance of L-DQN algorithm on \verb+cifar10+ is significantly better than other distributed algorithms including DAve-QN. L-DQN is competitive with DAve-QN and DANE on \verb+covtype+ as well.
\begin{figure}[ht]
    \centering
    \includegraphics[width=1\linewidth]{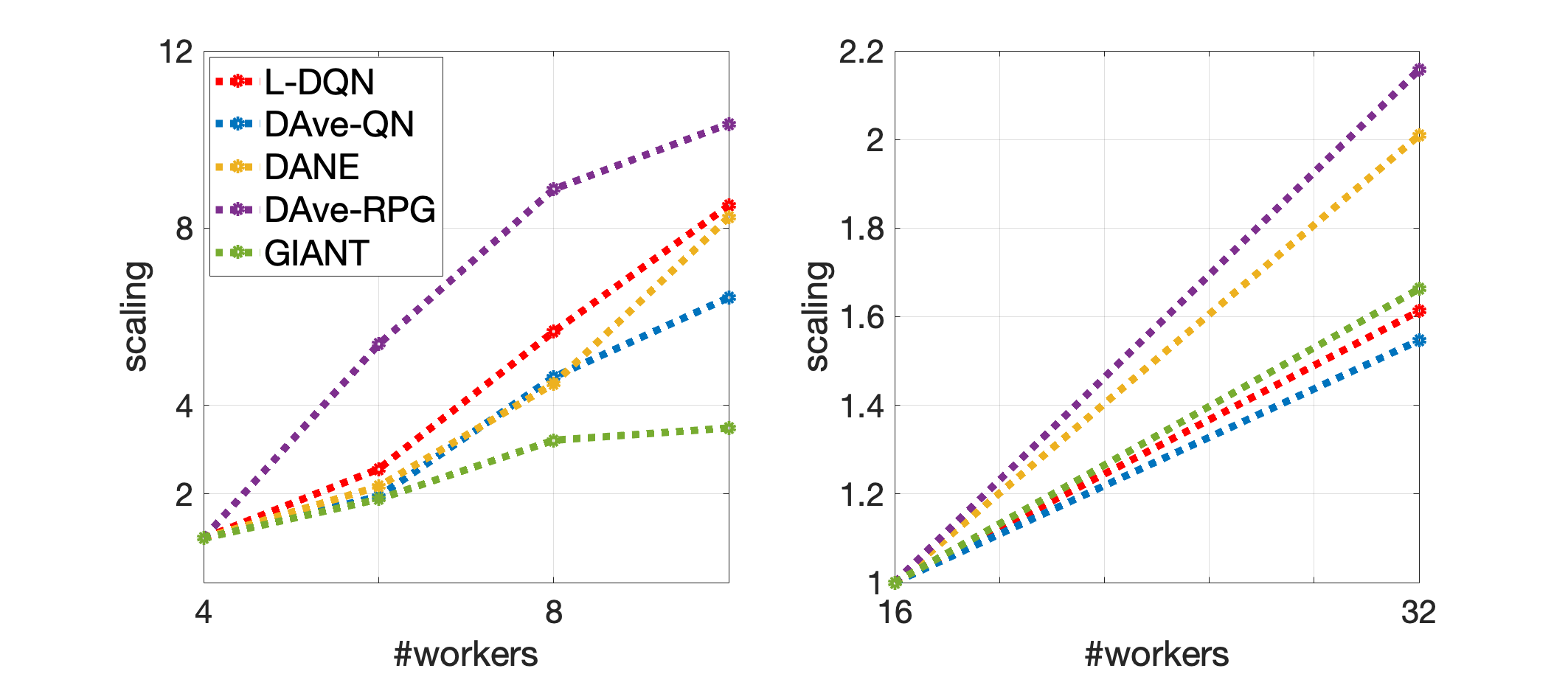}
    \caption{\small Scaling comparisons of the algorithms on \texttt{covtype}(left) and  \texttt{mnist8m}(right) datasets}
    \label{fig:Scaling Comparison}
\end{figure}

In Figure \ref{fig:Scaling Comparison}, we also compare the strong scaling of the distributed algorithms on different number of workers for \verb+mnist8m+ and \verb+covtype+. In particular, we look at the improvement in time to achieve the same suboptimality as we increase the number of workers. We see that L-DQN shows a nearly linear speedup and a slightly better scaling compared to DAve-QN. DAve-RPG scales better but considering the total runtime, it is slower than L-DQN.

In addition to suboptimality and scaling, we also compared the performance of these algorithms for different sparsity of the datasets. For the problem of interest (logistic regression), computing the gradient takes $O(nd)$ for dense and $O(n .\verb+nnz+)$ for sparse datasets where \verb+nnz+ is the number of non-zeros in the dataset. Therefore, L-DQN has  $O(n. \verb+nnz+ + md)$ while DAve-QN has a iteration complexity of $O(n d^2.\verb+nnz+ )$. Similarly, DAve-RPG has a complexity of $O(p . n . \verb+nnz+ +pd )$ where $p$ is number of passes on local data. We observe that L-DQN has a cheaper iteration complexity compared to DAve-QN while in case of very sparse datasets, DAve-RPG has a cheaper iteration complexity compared to L-DQN.This is illustrated over the dataset \verb+rcv1+ on the left panel of \autoref{fig:suboptimality-extra}. The dataset  \verb+rcv1+ is quite sparse with $\approx1\%$ non-zeros. We use the parameters $\lambda=0.01$, $m=10$,  $\eta=0.95$. For this dataset, DAve-QN fails as it requires more memory than the resources available. GIANT requires each worker to have $|S|>d$ where $|S|$ is the number of local data points on a worker. Hence, GIANT diverges with 16 workers. We observe that DAve-RPG converges faster than DANE and L-DQN because of its cheap iteration complexity. 

\begin{figure}[ht]
    \centering
    \includegraphics[width=0.6\linewidth]{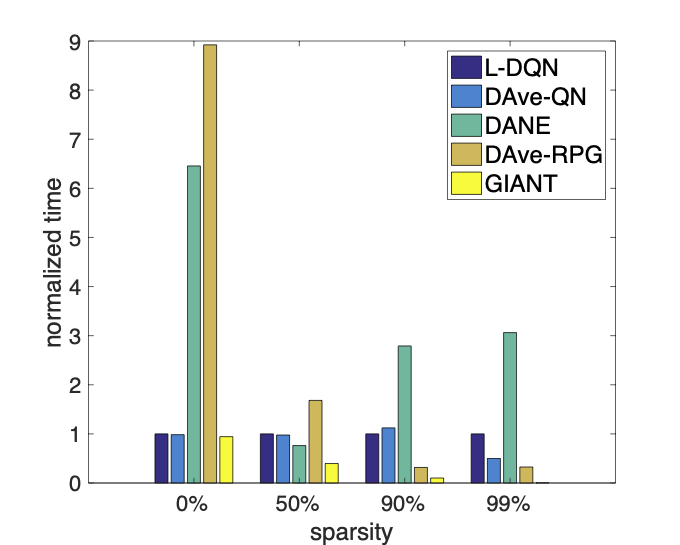}
    \caption{\small The effect of dataset sparsity on the performance of distributed optimization methods.}
    \label{fig:synth}
\end{figure}

In order to show the effect of sparsity on performance, we design a synthetic dataset based on a similar approach taken in \cite{shamir2014communication}. First we generate $N$ i.i.d input samples $x\sim \mathcal{N}(0, \Sigma)$ where $x\in\mathbb{R}^{2000}$ and the covariance matrix $\Sigma$ is diagonal with $\Sigma_{ii} = i^{-1.2}$. Then, we randomly choose some entries of all samples and make them zero to add sparsity. We set $z =  \langle\, x,w^*\rangle + \xi$, $\xi \sim \mathcal{N}(0, 0.09)$  and $w^*$ is the vector of all ones. Finally, labels $y \in\{0,1\}$ are generated based on the probabilities $p = S(z)$ where $S(z) = 1/(1+exp(-z))$ is the logistic function. The parameters $\lambda=0.01$ and $N=32000$ are chosen for the objective function and for this experiment we have the following $m=20$ and $\eta=0.9$. Time to the accuracy of $1e-4$ for all methods is measured and normalized based on L-DQN timing. The results are shown in \autoref{fig:synth}. DAve-RPG and DANE performs poorly for fully dense datasets (sparsity = 0\%), however, DAve-RPG and GIANT perform better compared to L-DQN as the dataset sparsity increases. We observe that when above \%90 of the data is sparse, DAve-RPG is the most efficient method; whereas for denser datasets GIANT and L-DQN are more efficient  on the synthetic data. 

\section{Conclusion}
We proposed the L-DQN method which is an asynchronous limited-memory BFGS method. We showed that under some assumptions, L-DQN admits linear convergence over epochs despite asynchronous computations. Our numerical experiments show that L-DQN can lead to significant performance improvements in practice in terms of both memory requirements and running time.


\bibliographystyle{unsrtnat}
\bibliography{main.bib}

\section*{Appendix}
\section{Proof of Theorem \ref{thm: L_DAveQN Convergence}}
Recall the definition of the average Hessian $\bar{G}_i^{t}= \int_0^1 \nabla^2f_i(x^*+\tau(z_i^t-x^*))(z_i^t-x^*)d\tau$ satisfies the equality $\nabla f_i(z_i^t)-\nabla f_i(x^*)=\bHi(z_i^t-x^*)$. Hence the equation \eqref{eq: Iterates} implies that iterates $x^t$ admit the bound
\vspace{-0.3cm}
\begin{small}
\begin{align}\label{Ineq: Iterate}
    \Vert x^t-x^*\Vert^2\leq \sum_{i=1}^{n}  \Vert \Gamma^t(\tB_i^t-\eta\bHi)\Vert^2_{2} \underset{i=1,..,n}{\max} \Vert z_i^t-x^*\Vert^2 
\end{align}
\end{small}
where $\Gamma^t={({\tB}^t)}^{-1}$ where ${\tB}^t$ is as in \eqref{update: Iterate_LDAveQN_Master} and $\Vert \cdot \Vert_2$ denotes the $2$-norm of a matrix. Notice that by its definition and from \eqref{Assump: Bounds}, it can be found that $\Gamma^t$ has the bounds $\frac{1}{n\lambda_u}I_d \preceq \Gamma^t \preceq \frac{1}{n\lambda_d}I_d$ and hence $(\Gamma^t)^{2}$ is positive definite. So the function $\Psi(A):= A(\Gamma^t)^{2}A$ defined from the set of symmetric positive-definite matrices $S^n$ to itself is a matrix convex function  \cite[see E.7.a]{MatrixInequalities}, that is for any $A,B\in S^n$ and $\alpha\in(0,1)$ following inequality holds 
\begin{equation}\label{ineq:matrix convexity} 
\Psi(\alpha A + (1-\alpha )B) \preceq \alpha \Psi(A)+(1-\alpha) \Psi(B).
\end{equation} 
In particular, if $L_A \preceq A \preceq U_A$ for some positive-definite matrices $L_A$ and $U_A$, by matrix convexity we have
\begin{small}
\begin{equation}\label{eq-mat-cvx-bounds}
    \sup_{L_A\preceq A \preceq U_A} \Psi(A) \preceq \max(\Psi(L_A), \Psi(U_A)).
\end{equation}
\end{small}
where the maximum on the right-hand side is in the sense of Loewner ordering, i.e. $\max\{A,B\}=A$ if $B\preceq A$ and equals to B otherwise. From the bounds \eqref{Assump: Bounds}, we have $
\left(1-\eta \frac{L}{\lambda_d}\right)\tB_i^t \preceq \tB_i^t-\eta\bHi \preceq \left(1-\eta\frac{\mu}{\lambda_u}\right) \tB_i^t,
$ for each  $i=1,..,n$. On the other hand, 
$
[\tB_i^t (\Gamma^t)^2 \tB_i^t]^{-1}= \Big(I_d+ \sum_{j\neq i}(\tB_i^t)^{-1}\tB_j \Big)\Big(I_d+ \sum_{j\neq i}\tB_j(\tB_i^t)^{-1} \Big)
$
together with \eqref{Assump: Bounds} imply that $\lambda_{\min}([\tB_i^t (\Gamma^t)^2 \tB_i^t]^{-1}) \geq  \left(1+(n-1)\frac{\lambda_d}{\lambda_u}\right)^2$, where $\mathrm{\lambda_{min}}$ is the smallest eigenvalue. This yields to
\begin{small}
\begin{align}
 \lambda_{\max}([\tB_i^t (\Gamma^t)^2 \tB_i^t])
 \leq  \frac{\lambda_u^2}{(\lambda_u+(n-1)\lambda_d)^2}.\label{ineq-lambda-max},
\end{align}
\end{small}
where $\mathrm{\lambda_{max}}$ denotes the largest eigenvalue. Applying \eqref{eq-mat-cvx-bounds} with $A=(\tB_i^t-\eta\bHi)$ with $L_A = \left(1-\eta \frac{L}{\lambda_d}\right)\tB_i^t$ and $U_A = \left(1-\eta\frac{\mu}{\lambda_u}\right) \tB_i^t$ and using \eqref{ineq-lambda-max}, we obtain
$\lambda_{\max}\Big((\tB_i^t-\eta\bHi)\Gamma^2 (\tB_i^t-\eta\bHi)\Big)\leq \frac{\lambda_u^2}{(\lambda_u+(n-1)\lambda_d)^2}\max \left\{ \left(1-\eta \frac{L}{\lambda_d}\right)^2,\left(1-\eta \frac{\mu}{\lambda_u}\right)^2 \right\}
$ for all $i=1,..,n$. 
Hence, 
\begin{footnotesize}
\begin{align}\label{ineq:rate}
&\sum_{i=1}^{n}\Vert \Gamma^t(\tB_i-\eta \bHi)\Vert^2_2\leq  \frac{n \tilde{\kappa}^2}{(\tilde{\kappa}+n-1)^2}\max\left\{ \left(1-\eta \frac{L}{\lambda_d} \right)^2,\left(1-\eta \frac{\mu}{\lambda_d} \right)^2 \right\}.
\end{align}
\end{footnotesize}
Choosing $\rho^2= \frac{n \tilde{\kappa}^2}{(\tilde{\kappa}+n-1)^2}\max\left\{ \left(1-\eta \frac{L}{\lambda_d} \right)^2,\left(1-\eta \frac{\mu}{\lambda_d} \right)^2 \right\}$ together with condition on $\eta$ imply that $\Vert x^t-x^*\Vert \leq \rho \max_{i=1,..,n}\Vert z_i^t-x^* \Vert$ where $\rho<1$. Next, we will prove convergence by induction on epoch times $E_{m}$. Notice that if $t\in [E_j,E_{j+1})$, it holds that $t-D_i^t \in [E_{j-1},t)$ for any $j\geq 1$, therefore  the inequality \eqref{ineq:rate} implies $\Vert x^t-x^*\Vert \leq \rho \max_{i=1,..,n}\Vert z_i^0-x^*\Vert\leq \rho\Vert x^0-x^*\Vert$ for $t\in [E_{0},E_1)$. Suppose for all $0\leq j \leq m$ the inequality $\Vert x^t-x^*\Vert \leq \rho^j \Vert x^0-x^*\Vert $ holds for $t\in [E_j,E_{j+1})$, then \eqref{Ineq: Iterate} and \eqref{ineq:rate} imply $ \Vert x^{t}-x^* \Vert \leq \rho \max_{i=1,..,n}\Vert z_i^{t-D_i^t}-x^*\Vert \leq  \rho^{m}\Vert x^0-x^*\Vert.$ This completes the proof.

\end{document}